\documentclass[12pt]{article}
\usepackage{cancel}
\usepackage{tabularx}
\usepackage{tikz}
\usepackage{relsize}
mm \textheight=8.9in

\usepackage{amssymb}
\usepackage{amsmath}
\usepackage{amsthm}
\usepackage{color}

\usepackage{tikz}
\usetikzlibrary{patterns,snakes}
\usepackage{pict2e}

\usepackage{tikz-cd}

\usepackage{array,multirow}

\newcommand{\sdim}{\operatorname{sdim}}
\newcommand{\op}{\mathrm{op}}

\newcommand{\br}[3]{{$#1$}$\lower4pt\hbox{$\tp\atop\raise4pt \hbox{$\scriptscriptstyle{#2}$}$} ${$#3$}}
\newcommand{\tw}[3]{{$#1$}${\,\scriptscriptstyle {#2}}\atop\raise9pt\hbox{$\scriptstyle\tp$} ${$#3$}}
\newcommand{\ttps}[2]{{#1}\raise5pt\hbox{$\lower12pt\hbox{$\scriptstyle\tp$}\atop \lower0pt\hbox{$\tilde\;$}$}\raise4.5pt\hbox{${\scriptstyle{#2}}$}}
\newcommand{\st}[1]{\mbox{${\,\scriptscriptstyle {#1}}\atop\raise5.5pt\hbox{$*$}$}}

\newcommand{\rd}[1]{\mbox{${\,\scriptscriptstyle {#1}}\atop\raise5.5pt\hbox{$\bullet$}$}}
\newcommand{\rt}[1]{\otimes_\chi}
\newcommand{\lt}[1]{\mbox{${\,\scriptscriptstyle {#1}}\atop\raise5.5pt\hbox{$\ltimes$}$}}
\newcommand{\btr}{\raise1.2pt\hbox{$\scriptstyle\blacktriangleright$}\hspace{2pt}}
\newcommand{\btl}{\raise1.2pt\hbox{$\scriptstyle\blacktriangleleft$}\hspace{2pt}}

\newcommand{\lcr}{\raise1.0pt \hbox{${\scriptstyle\rightharpoonup}$}}
\newcommand{\rcr}{\raise1.0pt \hbox{${\scriptstyle\leftharpoonup}$}}

\newcommand{\ttp}{{\lower12pt\hbox{$\tp$}\atop \hbox{$\tilde\;$}}}

\newcommand{\id}{\mathrm{id}}

\newcommand{\Ru}{\mathcal{R}}

\newcommand{\C}{\mathbb{C}}

\newcommand{\tp}{\otimes}

\newcommand{\End}{\mathrm{End}}

\newcommand{\g}{\mathfrak{g}}

\newcommand{\h}{\mathfrak{h}}

\newcommand{\s}{\mathfrak{s}}

\renewcommand{\o}{\mathfrak{o}}

\newcommand{\p}{\mathfrak{p}}
\renewcommand{\l}{\mathfrak{l}}

\newcommand{\al}{\alpha}

\newcommand{\be}{\begin{eqnarray}}
\newcommand{\ee}{\end{eqnarray}}

\newtheorem{thm}{Theorem}[section]
\newtheorem{propn}[thm]{Proposition}
\newtheorem{lemma}[thm]{Lemma}

\newcount\prg

\newcommand{\parag}{\advance\prg by1 {\noindent\bf\thesection.\the\prg\hspace{6pt}}}

\begin{document}

\title{Orthosymplectic $R$-matrix from Hasse diagrams}
\author{D. Algethami${}^{\dag}$, N.~Kryazhevskikh${}^{\ddag}$, and A. Mudrov${}^{\ddag}$
 \vspace{10pt}\\
\small ${\dag}$ Department of Mathematics, College of Science,\\
\small
 University of Bisha, P.O. Box 551, Bisha 61922, Saudi Arabia,
\vspace{10pt}\\
\small
${\ddag}$
\small
Centre for Pure Mathematics\\
\small
 Moscow Institute of Physics and Technology,\\
\small 9 Institutskiy per., Dolgoprudny, Moscow Region,
141701, Russia,
\vspace{10pt}\\
\small
     e-mail: dalgethami@ub.edu.sa, krjazhevskikh.nv@phystech.edu, mudrov.ai@mipt.ru
}

\date{}
\maketitle

\begin{abstract}
We calculate R-matrix in the basic representation   for the ortho-symplectic quantum supergroups in  arbitrary polarizations,
by a method of Hasse diagrams.
\end{abstract}

{\small \underline{2010 AMS Subject Classification}: 17B37, 17A70.\\[4pt]
}

\newpage
\section{Introduction}
The Yang--Baxter equation \cite{B} is a key concept in a number of mathematical theories including quantum groups \cite{D1,FRT}, integrable models \cite{FT}, and knot invariants \cite{RT}.
Its solutions called R-matrices deliver symmetry morphisms in relative quasi-tensor categories. The quantum groups
were    specifically designed to accommodate the  collection of intertwiners  in a universal R-matrix $\Ru$ \cite{D1},
which serves as a source of solutions in particular representations. Of special interest are R-matrices in fundamental
representations, including those of minimal dimension and the adjoint representation.
Explicit factorization of $\Ru$ over positive roots obtained by several authors \cite{KR,KT,LS, Y2} facilitates such calculations
 in principle. However there remains a problem of an efficient presentation of the result as expansion over a natural matrix basis.

An alternative way of constructing R-matrices is via a recursion over Hasse diagrams associated with modules diagonalizable over the Cartan
subalgebras with bounded weights.
This approach gives solutions for  modules of highest weight and for finite dimensional modules in particular \cite{KMS}.
In a relatively simple situation as the natural modules of the basic matrix supergroups it gives a straightforward and easily verifiable
computational algorithm.
One has to bear in mind that simple Lie superalgebras have a variety of generally non-isomorphic Borel subalgebras determined
by a parity sequence on basic weights. Different Borel subalgebras give rise to different fundamental R-matrices.
Until recently, such an R-matrix had been known only for the distinguished polarization (with one odd simple root)
\cite{DGL}.
Hong and Tsymbaliuk found it for an arbitrary polarization using the universal R-matrix approach \cite{HT}.
Here we derive these results via the method of Hasse diagrams.

In the next section, we recall the basic representation of orthosymplectic quantum groups, with the same conventions as in \cite{HT}.
We then describe the Hasse diagrams and the recursion algorithm and calculate the orthosymplectic R-matrices.
In the last section, we compare the result with \cite{HT} and derive a cubic relation for the braid matrix, including its minimal polynomial at zero superdimension.

\section{The natural representation}

Let $V=\C^N$, $N>2$, be a vector superspace with a homogeneous basis
$v_1,\ldots,v_N$. Set
\[
 n=\lfloor N/2\rfloor,\qquad i'=N+1-i,\qquad
 p_i=|v_i|\in\{0,1\},\qquad s_i=(-1)^{p_i}.
\]
We assume that $p_{i'}=p_i$ and, if $N=2n+1$, that $p_{n+1}=0$.
The matrix units $e_{ij}$ are defined by their action  $e_{ij}v_k=\delta_{jk}v_i$
and have parity $p_i+p_j$. All parity expressions are understood modulo $2$.
Fix the signs
\begin{equation}\label{eq:theta}
 \vartheta_i=1\quad(i\le\lceil N/2\rceil),\qquad
 \vartheta_{i'}=s_i\quad(i\le n).
\end{equation}
Then $\vartheta_i\vartheta_{i'}=s_i$, and the bilinear form
\[
 \langle v_i,v_j\rangle=\vartheta_j\delta_{i,j'}
\]
is even, nondegenerate, and supersymmetric:
$\langle v_i,v_j\rangle=(-1)^{p_ip_j}\langle v_j,v_i\rangle$.
The Lie superalgebra preserving this form is denoted by
$\g=\mathfrak{osp}(V)$. For homogeneous $A\in\g$ and $v,w\in V$,
the invariance condition reads
$\langle Av,w\rangle+(-1)^{|A||v|}\langle v,Aw\rangle=0$.

Choose the diagonal Cartan subalgebra $\h\subset\g$ and the upper
triangular Borel subalgebra. Write $\zeta_i$ for the weight of $v_i$
and normalize the inner product by
\begin{equation}\label{eq:weights}
 \zeta_{i'}=-\zeta_i,\qquad
 (\zeta_i,\zeta_j)=s_i(\delta_{ij}-\delta_{i,j'}).
\end{equation}
Thus $\zeta_{n+1}=0$ for odd $N$, and
$(\zeta_i,\zeta_i)=s_i$ for $i\ne i'$. The simple roots are
\begin{equation}\label{eq:roots}
 \alpha_i=\zeta_i-\zeta_{i+1}\quad(i<n),\qquad
 \alpha_n=
 \begin{cases}
 \zeta_n,&N=2n+1\quad(B),\\
 2\zeta_n,&N=2n,\quad(C),\\
 \zeta_{n-1}+\zeta_n,&N=2n,\quad(D).
 \end{cases}
\end{equation}
The letters $B,C,D$ indicate the terminal type of the diagram.
For even $N$, the terminal type is $C$ if $p_n=1$ and $D$ if $p_n=0$.
The parity of $\alpha_i$ is $p_i+p_{i+1}$ for $i<n$;
that of $\alpha_n$ is $p_n$, $0$, or $p_{n-1}+p_n$ in types $B,C,D$,
respectively.

We work over $\C(q^{1/2})$, with $q$ indeterminate, and write
$\bar q=q^{-1}$ and $2_q=q+\bar q$.
The quantum superalgebra $U_q(\g)$ is generated by
$e_i,f_i,q^{\pm h_i/2}$, where $h_i$ corresponds to $\alpha_i$
under the inner product. Our conventions are
\begin{equation}\label{eq:relations}
 \begin{gathered}
 q^{h_i/2}e_jq^{-h_i/2}=q^{(\alpha_i,\alpha_j)/2}e_j,\qquad
 q^{h_i/2}f_jq^{-h_i/2}=q^{-(\alpha_i,\alpha_j)/2}f_j,\\
 \{e_i,f_j\}=\delta_{ij}\frac{q^{h_i}-q^{-h_i}}{q-\bar q}.
 \end{gathered}
\end{equation}
Here $\{x,y\}=xy-(-1)^{|x||y|}yx$ denotes the supercommutator.
The generators $e_i,f_i$ have the parity of $\alpha_i$.
The elements $q^{\pm h_i/2}$ are even, invertible, and pairwise
commuting. The defining relations also include the standard and
additional quantum Serre relations \cite[Section~2.4]{HT}.

We use the natural representation
$\pi:U_q(\g)\to\End(V)$ of \cite[Proposition~3.1]{HT},
with the signs~\eqref{eq:theta}. The Cartan generators act by
\begin{equation}\label{eq:cartan}
 \pi(q^{\pm h_i/2})=\sum_j q^{\pm(\alpha_i,\zeta_j)/2}e_{jj}.
\end{equation}
For $i<n$, the remaining generators are
\begin{equation}\label{eq:chain}
 \begin{split}
 \pi(e_i)&=e_{i,i+1}
        -(-1)^{p_i(p_i+p_{i+1})}e_{(i+1)',i'},\\
 \pi(f_i)&=s_i\bigl(e_{i+1,i}
        -(-1)^{p_{i+1}(p_i+p_{i+1})}e_{i',(i+1)'}\bigr).
 \end{split}
\end{equation}
At the terminal root, they are given by
\begin{equation}\label{eq:terminal}
 \begin{array}{c|c|c}
 &\pi(e_n)&\pi(f_n)\\ \hline
 B&e_{n,n+1}-s_ne_{n+1,n'}&s_n(e_{n+1,n}-e_{n',n+1})\\[3pt]
 C&e_{n,n'}&-2_qe_{n',n}\\[3pt]
 D&e_{n-1,n'}-s_{n-1}e_{n,(n-1)'}&
        s_{n-1}(e_{n',n-1}-e_{(n-1)',n}).
 \end{array}
\end{equation}
In particular, $v_1$ is a highest vector of weight $\zeta_1$.

We check the normalizations used below. For $i<n$, multiplication gives
\[
 \begin{split}
 \pi(e_i)\pi(f_i)&=s_ie_{ii}+s_{i+1}e_{(i+1)',(i+1)'},\\
 \pi(f_i)\pi(e_i)&=s_ie_{i+1,i+1}+s_{i+1}e_{i'i'}.
 \end{split}
\]
Since $(-1)^{p_i+p_{i+1}}=s_is_{i+1}$, it follows that
\[
 \{\pi(e_i),\pi(f_i)\}
 =s_i(e_{ii}-e_{i'i'})
       -s_{i+1}(e_{i+1,i+1}-e_{(i+1)',(i+1)'}).
\]
By~\eqref{eq:weights}, its $j$th diagonal entry is
\[
 \begin{split}
 &s_i(\delta_{ij}-\delta_{i'j})
       -s_{i+1}(\delta_{i+1,j}-\delta_{(i+1)'j})\\
 &\qquad=(\zeta_i,\zeta_j)-(\zeta_{i+1},\zeta_j)
          =(\alpha_i,\zeta_j).
 \end{split}
\]
These entries lie in $\{0,\pm1\}$, so the supercommutator equals
$\bigl(\pi(q^{h_i})-\pi(q^{-h_i})\bigr)/(q-\bar q)$.
At the ends of types $B,D$, the corresponding products give
\[
 \{\pi(e_n),\pi(f_n)\}=
 \begin{cases}
 s_n(e_{nn}-e_{n'n'}),&B,\\
 s_{n-1}(e_{n-1,n-1}-e_{(n-1)',(n-1)'})
                   +e_{nn}-e_{n'n'},&D.
 \end{cases}
\]
Their diagonal entries are again $(\alpha_n,\zeta_j)\in\{0,\pm1\}$.
For the long simple root of type $C$, we have $(\alpha_n,\zeta_n)=-2$ and
\[
 \{\pi(e_n),\pi(f_n)\}
 =-2_q(e_{nn}-e_{n'n'})
 =\frac{\pi(q^{h_n})-\pi(q^{-h_n})}{q-\bar q}.
\]
This accounts for the factor $2_q$ in~\eqref{eq:terminal}.
Each matrix unit in $\pi(e_i)$ has weight $\alpha_i$, and each
matrix unit in $\pi(f_i)$ has weight $-\alpha_i$, which also verifies
the conjugation relations in~\eqref{eq:relations}.
The Serre relations for these matrices are established in
\cite[Proposition~3.1]{HT}.

All tensor products are graded. For homogeneous elements,
\begin{equation}\label{eq:graded}
 (A\otimes B)(C\otimes D)=(-1)^{|B||C|}AC\otimes BD.
\end{equation}
In particular,
\[
 (e_{ij}\otimes e_{kl})(e_{ab}\otimes e_{cd})
 =(-1)^{(p_k+p_l)(p_a+p_b)}
                      \delta_{ja}\delta_{lc}\,e_{ib}\otimes e_{kd}.
\]
The graded permutation is $P=\sum_{i,j}s_je_{ij}\otimes e_{ji}$.
It satisfies $P^2=\id$ and
$P(A\otimes B)P=(-1)^{|A||B|}B\otimes A$.
We use the symmetric coproduct
\begin{equation}\label{eq:coproduct}
 \begin{split}
 \Delta(e_i)&=q^{h_i/2}\otimes e_i+e_i\otimes q^{-h_i/2},\\
 \Delta(f_i)&=q^{h_i/2}\otimes f_i+f_i\otimes q^{-h_i/2},\\
 \Delta(q^{\pm h_i/2})&=q^{\pm h_i/2}\otimes q^{\pm h_i/2},
 \end{split}
\end{equation}
as in \cite[(2.23)]{HT}. In formulas for operators on $V\otimes V$,
we abbreviate $(\pi\otimes\pi)\Delta(x)$ to $\Delta(x)$ and set
$\Delta^{\op}(x)=P\Delta(x)P$.
\section{Recursion along the Hasse diagram}

We compute the image of the normalized universal $R$-matrix in
$\End(V)\otimes\End(V)$. It satisfies the intertwining equation

\begin{equation}\label{eq:intertwining}
 R\Delta(x)=\Delta^{\op}(x)R,\qquad x\in U_q(\g).
\end{equation}
With the normalization of \cite[Section~4.1]{HT}, its diagonal part is
\[
 \sum_{i,j}q^{-(\zeta_i,\zeta_j)}e_{ii}\otimes e_{jj},
\]
and its off-diagonal terms have weights that are nonzero sums of
positive roots in the first tensor factor and opposite weights in
the second. Write
\begin{equation}\label{eq:blocks}
 R=\sum_{i,j}e_{ij}\otimes L_{ij},\qquad L_{ij}\in\End(V).
\end{equation}
The coefficient $L_{ij}$ has parity $p_i+p_j$ and weight
$\zeta_j-\zeta_i$. Triangularity gives $L_{ij}=0$ for $i>j$, while
the first column is
\begin{equation}\label{eq:first-column}
 L_{i1}=\delta_{i1}\sum_j q^{-(\zeta_1,\zeta_j)}e_{jj}.
\end{equation}

The following Hasse diagrams describe the natural representation
restricted to the subalgebra $U_q(\g_-)\subset U_q(\g)$ generated by
$f_1,\ldots,f_n$. They correspond to the terminal types $B,C,D$,
respectively; vertex parities are omitted.
 \begin{center}
\begin{picture}(390,45)
\put(8,0){${\scriptstyle 1'}$}
\put(10,15){\circle{3}}
\put(45,15){\vector(-1,0){30}}
\put(50,15){\circle{3}}
\put(100,15){$\ldots$}
\put(85,15){\vector(-1,0){30}}
\put(47,0){${\scriptstyle 2'}$}
\put(160,15){\circle{3}}
\put(155,15){\vector(-1,0){30}}
\put(158,0){${\scriptstyle n'}$}
\put(195,15){\vector(-1,0){30}}
\put(200,15){\circle{3}}
\put(194,0){${\scriptstyle n+1}$}
\put(235,15){\vector(-1,0){30}}
\put(240,15){\circle{3}}
\put(239,0){${\scriptstyle n}$}
\put(275,15){\vector(-1,0){30}}
\put(345,15){\vector(-1,0){30}}
\put(350,15){\circle{3}}
\put(348,0){${\scriptstyle 2}$}
\put(290,15){$\ldots$}
\put(385,15){\vector(-1,0){30}}
\put(390,15){\circle{3}}
\put(388,0){${\scriptstyle 1}$}
\put(25,20){${\scriptstyle\al_1}$}
\put(65,20){${\scriptstyle\al_2}$}
\put(135,20){${\scriptstyle\al_{n-1}}$}
\put(175,20){${\scriptstyle\al_{n}}$}
\put(215,20){${\scriptstyle\al_{n}}$}
\put(255,20){${\scriptstyle\al_{n-1}}$}
\put(325,20){${\scriptstyle\al_{2}}$}
\put(365,20){${\scriptstyle\al_{1}}$}

\end{picture}
\end{center}

\begin{center}
\begin{picture}(390,45)
\put(8,0){${\scriptstyle 1'}$}
\put(10,15){\circle{3}}
\put(45,15){\vector(-1,0){30}}
\put(50,15){\circle{3}}
\put(100,15){$\ldots$}
\put(85,15){\vector(-1,0){30}}
\put(48,0){${\scriptstyle 2'}$}
\put(160,15){\circle{3}}
\put(155,15){\vector(-1,0){30}}
\put(158,0){${\scriptstyle n'}$}
\put(195,15){\vector(-1,0){30}}
\put(200,15){\circle{3}}
\put(198,0){${\scriptstyle n}$}
\put(235,15){\vector(-1,0){30}}
\put(305,15){\vector(-1,0){30}}
\put(310,15){\circle{3}}
\put(308,0){${\scriptstyle 2}$}
\put(250,15){$\ldots$}
\put(345,15){\vector(-1,0){30}}
\put(350,15){\circle{3}}
\put(348,0){${\scriptstyle 1}$}
\put(25,20){${\scriptstyle\al_1}$}
\put(65,20){${\scriptstyle\al_2}$}
\put(135,20){${\scriptstyle\al_{n-1}}$}
\put(175,20){${\scriptstyle\al_{n}}$}
\put(215,20){${\scriptstyle\al_{n-1}}$}
\put(285,20){${\scriptstyle\al_{2}}$}
\put(325,20){${\scriptstyle\al_{1}}$}

\end{picture}
\end{center}

\begin{center}
\begin{picture}(440,50)
\put(8,0){${\scriptstyle 1'}$}
\put(10,15){\circle{3}}
\put(50,15){\circle{3}}
\put(45,15){\vector(-1,0){30}}
\put(85,15){\vector(-1,0){30}}
\put(100,15){$\ldots$}
\put(155,15){\vector(-1,0){30}}
\put(48,0){${\scriptstyle 2'}$}
\put(160,15){\circle{3}}
\put(149,0){${\scriptstyle (n-1)'}$}
\put(195,15){\vector(-1,0){30}}
\put(200,15){\circle{3}}
\put(197,0){${\scriptstyle n'}$}
\qbezier(163,22)(200,48)(237,22) \put(162,21.5){\vector(-2,-1){2}}
\qbezier(203,22)(240,48)(277,22) \put(202,21.5){\vector(-2,-1){2}}
\put(240,15){\circle{3}}
\put(238,0){${\scriptstyle n}$}
\put(275,15){\vector(-1,0){30}}
\put(280,15){\circle{3}}
\put(273,0){${\scriptstyle n-1}$}
\put(330,15){$\ldots$}
\put(315,15){\vector(-1,0){30}}
\put(385,15){\vector(-1,0){30}}
\put(390,15){\circle{3}}
\put(388,0){${\scriptstyle 2}$}
\put(425,15){\vector(-1,0){30}}
\put(430,15){\circle{3}}
\put(428,0){${\scriptstyle 1}$}
\put(25,20){${\scriptstyle\al_1}$}
\put(65,20){${\scriptstyle\al_2}$}
\put(135,20){${\scriptstyle\al_{n-2}}$}
\put(172,20){${\scriptstyle\al_{n-1}}$}
\put(195,38){${\scriptstyle\al_{n}}$}
\put(252,20){${\scriptstyle\al_{n-1}}$}
\put(290,20){${\scriptstyle\al_{n-2}}$}
\put(235,38){${\scriptstyle\al_{n}}$}
\put(365,20){${\scriptstyle\al_{2}}$}
\put(405,20){${\scriptstyle\al_{1}}$}
\end{picture}
\end{center}
The vertices represent the homogeneous weight basis elements;
an arrow $b\xleftarrow{\alpha_k}a$ means $\pi(f_k)_{ba}\ne0$.
When $V$ has both even and odd basis vectors, the second and third
diagrams arise from different polarizations of the same Lie superalgebra.

Every vertex is reachable from $1$. For a fixed simple root,
at most one arrow leaves any given vertex.
For an arrow $b\xleftarrow{\alpha_k}a$, the parity of $f_k$ is
$p_a+p_b$. Comparing coefficients of $e_{ia}$ in the first tensor
factor of~\eqref{eq:intertwining} with $x=f_k$ gives
\begin{equation}\label{eq:hasse}
 \begin{split}
 &(-1)^{(p_i+p_b)(p_a+p_b)}
       \pi(f_k)_{ba}L_{ib}\pi(q^{-h_k/2})\\
 &\quad=(-1)^{(p_i+p_a)(p_a+p_b)}
       q^{-(\alpha_k,\zeta_i)/2}\pi(f_k)L_{ia}\\
 &\qquad\quad-q^{(\alpha_k,\zeta_a)/2}L_{ia}\pi(f_k)
       +\sum_j\pi(f_k)_{ij}\pi(q^{h_k/2})L_{ja}.
 \end{split}
\end{equation}
All coefficients on the right belong to column $a$.
Since $\pi(f_k)_{ba}\ne0$ and $\pi(q^{-h_k/2})$ is invertible,
this determines column $b$.

\begin{propn}\label{prop:unique}
The intertwining equations~\eqref{eq:intertwining} for the simple
lowering generators and the initial column~\eqref{eq:first-column}
determine $R$ uniquely.
\end{propn}
\begin{proof}
Apply~\eqref{eq:hasse} successively along a path from $1$ to each
vertex. Existence of the normalized universal $R$-matrix ensures
compatibility, including the two paths through the fork of type $D$.
\end{proof}

For distinct nonopposite indices $a<b$, weight preservation and
triangularity restrict the corresponding block of the R-matrix to
\begin{equation}\label{eq:exchange}
 e_{aa}\otimes e_{bb}+e_{bb}\otimes e_{aa}
                         +x_{ab}e_{ab}\otimes e_{ba}.
\end{equation}
Here $x_{ab}$ are exchange coefficients to be found.
The coefficient of $e_{ii}\otimes e_{ii}$ is $q^{-s_i}$ for $i\ne i'$.

\begin{lemma}\label{lem:exchange}
For $a<b$ and $b\ne a'$, the coefficient in~\eqref{eq:exchange} is
$x_{ab}=-(q-\bar q)s_b$.
\end{lemma}
\begin{proof}
For a simple arrow $b\xleftarrow{\alpha_k}a$ with $\zeta_a\ne0$,
comparison of the coefficients of $e_{aa}\otimes e_{ba}$ gives
\[
 q^{(\zeta_a,\zeta_a)/2}
 +(-1)^{p_a+p_b}q^{-(\zeta_a,\zeta_a)/2}x_{ab}
 =q^{-3(\zeta_a,\zeta_a)/2}.
\]
Hence $(-1)^{p_a+p_b}x_{ab}=-s_a(q-\bar q)$, as required.
For the arrow $n'\leftarrow n+1$ of type $B$, comparison at
$e_{bb}\otimes e_{ba}$ in (\ref{eq:intertwining}) gives
$x_{ab}=q^{-(\zeta_b,\zeta_b)}-q^{(\zeta_b,\zeta_b)}=-s_b(q-\bar q)$.

If $d$ lies outside the row and column support of $\pi(f_k)$,
comparison at $e_{da}\otimes e_{bd}$ and $e_{bd}\otimes e_{da}$ gives
\begin{equation}\label{eq:exchange-transport}
 x_{db}=(-1)^{p_a+p_b}x_{da}\quad(b>a>d),\qquad
 x_{bd}=x_{ad}\quad(d>b>a).
\end{equation}
Thus $s_jx_{ij}$ is preserved by either admissible shift of the pair
$j>i$. In the chains of types $B,C$, a pair with $j-i>1$ can be
shortened by replacement $i\to i+1$ or $j\to j-1$: the support condition fails
for the first shift only if $i+j=N$, and for the second only if
$i+j=N+2$. The opposite pair $i+j=N+1$ is excluded.
In type $D$, the only missing adjacent edge is $n'\leftarrow n$.
The arrows $n+2\leftarrow n$ and $n'\leftarrow n-1$ give the
terminal cases. Pairs $i=n<j$ and $i<j=n'$ reduce to these by
shifting the other endpoint. This reduces every nonopposite pair
to a simple arrow.
\end{proof}
\section{Contraction coefficients}

The remaining terms are $e_{i'j'}\otimes e_{ij}$, with both row and
column of total weight zero. Extending the exchange sum to opposite
indices, write
\begin{equation}\label{eq:R-ansatz}
 \begin{split}
 R={}&\sum_{i,j}q^{-(\zeta_i,\zeta_j)}e_{ii}\otimes e_{jj}
          -(q-\bar q)\sum_{i>j}s_i e_{ji}\otimes e_{ij}\\
 &+(q-\bar q)\sum_{i>j}(-1)^{p_j(p_i+p_j)}
                              c_{ij}e_{i'j'}\otimes e_{ij},
 \end{split}
\end{equation}
where the contraction coefficients $c_{ij}$ are to be determined.
For $i=j'$, the exchange and contraction terms coincide and their
coefficients add. Set $c_{jj}=s_j$ to include the initial conditions
in the recursion.

The paired arrows in~\eqref{eq:chain}--\eqref{eq:terminal} have the form
\begin{equation}\label{eq:paired-generator}
 \pi(f_k)=s_a\bigl(e_{ba}-(-1)^{p_b(p_a+p_b)}e_{a'b'}\bigr),
 \qquad \alpha_k=\zeta_a-\zeta_b,
\end{equation}
where $(b,a) =(k+1,k)$ for $k<n$, $(n+1,n)$ at the end of type $B$,
or $(n',n-1)$ at the end of type $D$.
For the first two cases, comparison at $e_{a'j'}\otimes e_{bj}$
in~\eqref{eq:intertwining} yields
\begin{equation}\label{eq:c-forward}
 c_{bj}=(-1)^{p_a+p_b}
          q^{((\zeta_a,\zeta_a)+(\zeta_b,\zeta_b))/2}c_{aj},
 \qquad j\le a.
\end{equation}
Comparison at $e_{bj'}\otimes e_{a'j}$ for the reflected arrows gives
\begin{equation}\label{eq:c-reflected}
 c_{a'j}=q^{((\zeta_a,\zeta_a)+(\zeta_b,\zeta_b))/2}c_{b'j},
 \qquad j\le b'.
\end{equation}
At the initial steps, the diagonal terms give
\[
 c_{ba}=c_{a'b'}=s_bq^{((\zeta_a,\zeta_a)+(\zeta_b,\zeta_b))/2}.
\]
Thus~\eqref{eq:c-forward} and~\eqref{eq:c-reflected} also hold
for $j=a$ and $j=b'$, respectively.
In type $B$, the central step $b=b'=n+1$ gives
$c_{a'b}=q^{(\zeta_a,\zeta_a)/2}$, also as in~\eqref{eq:c-reflected}.

In type $C$, comparison at $e_{n',j'}\otimes e_{n',j}$, after
cancelling the factor $-2_q$, gives
\[
 \bar q c_{nj}+q c_{n'j}=0\quad(j<n),\qquad
 \bar q^2+q(q-\bar q)(c_{n'n}+1)=q^2.
\]
Since $c_{nn}=-1$, these are equivalent to
\begin{equation}\label{eq:c-C}
 c_{n'j}=-q^{-2}c_{nj}\qquad(j\le n).
\end{equation}
In type $D$, comparison at $e_{(n-1)',j'}\otimes e_{nj}$ for
$f_{n-1}$ and at $e_{(n-1)',j'}\otimes e_{n'j}$ for $f_n$ gives
identical equations for $j<n$ after replacing $c_{nj}$ by $c_{n'j}$.
For $j=n$, the latter equation is
$q^{-1/2}(q-\bar q)(c_{n'n}-1)=0$. Thus
\begin{equation}\label{eq:c-D}
 c_{n'j}=c_{nj}\qquad(j\le n).
\end{equation}
Relations~\eqref{eq:c-forward}--\eqref{eq:c-D}, starting from
$c_{jj}=s_j$, determine all the coefficients.

To express their products, introduce the Weyl vector
\[
 \rho=\frac12\sum_{\substack{\alpha>0\\\alpha\text{ even}}}\alpha
        -\frac12\sum_{\substack{\alpha>0\\\alpha\text{ odd}}}\alpha.
\]
The sums run over the positive roots of the chosen Borel subalgebra.
The identities $(\rho,\alpha_i)=(\alpha_i,\alpha_i)/2$ give
$(\rho,\zeta_i-\zeta_{i+1})=(s_i+s_{i+1})/2$ for $i<n$.
The terminal value $(\rho,\zeta_n)$ is $s_n/2$, $-1$, or $0$
in types $B,C,D$, respectively. Consequently,
\begin{equation}\label{eq:weyl}
 \begin{split}
 (\rho,\zeta_i)
   &=\sum_{k=i+1}^n s_k+
     \begin{cases}
       s_i/2,&N=2n+1,\\
       (s_i-1)/2,&N=2n,
     \end{cases}\qquad i\le n,\\
 (\rho,\zeta_{i'})&=-(\rho,\zeta_i).
 \end{split}
\end{equation}
For odd $N$, we also have $(\rho,\zeta_{n+1})=0$.

\begin{propn}\label{prop:coefficients}
The solution of the recursion is
\begin{equation}\label{eq:c-closed}
     c_{ij}=s_i\vartheta_i\vartheta_jq^{(\rho,\zeta_j-\zeta_i)}
 \qquad(i\ge j).
\end{equation}
\end{propn}
\begin{proof}
The initial value is $s_j$. Along a direct arrow $b\leftarrow a$,
$\vartheta_a=\vartheta_b=1$ and $s_b/s_a=(-1)^{p_a+p_b}$;
by~\eqref{eq:weyl}, the ratio in~\eqref{eq:c-closed} is the multiplier
in~\eqref{eq:c-forward}. Along the reflected arrow,
$s_{a'}\vartheta_{a'}=s_{b'}\vartheta_{b'}=1$ and
$\zeta_{b'}-\zeta_{a'}=\zeta_a-\zeta_b$, giving~\eqref{eq:c-reflected}.
This includes the central vertex of type $B$, where
$s_{n+1}=\vartheta_{n+1}=1$.
The ratio from $n$ to $n'$ is
\[
 \frac{s_{n'}\vartheta_{n'}}{s_n\vartheta_n}
       q^{(\rho,\zeta_n-\zeta_{n'})}
 =\begin{cases}-q^{-2},&C,\\1,&D,\end{cases}
\]
which verifies the remaining terminal conditions.
The formula therefore follows by induction along the recursion.
\end{proof}
\section{The explicit matrix and operator identities}

Set $I=\id_{V\otimes V}$. Combining the coefficients gives the
following formula.

\begin{thm}\label{thm:R}
In the natural representation~\eqref{eq:cartan}--\eqref{eq:terminal},
with coproduct~\eqref{eq:coproduct}, the normalized $R$-matrix is
\begin{equation}\label{eq:R}
 \begin{split}
 R={}&I+(q^{-1/2}-q^{1/2})\sum_{i=1}^N
       s_ie_{ii}\otimes
       \bigl(q^{-(\zeta_i,\zeta_i)/2}e_{ii}
                      -q^{(\zeta_i,\zeta_i)/2}e_{i'i'}\bigr)\\
 &-(q-\bar q)\sum_{i<j}s_je_{ij}\otimes
       \bigl(e_{ji}-(-1)^{p_j(p_i+p_j)}
          \vartheta_i\vartheta_jq^{(\rho,\zeta_i-\zeta_j)}e_{i'j'}\bigr).
 \end{split}
\end{equation}
\end{thm}
\begin{proof}
The first line has coefficients $q^{-(\zeta_i,\zeta_j)}$ at
$e_{ii}\otimes e_{jj}$, since for $i\ne i'$,
\[
 1+s_i(q^{-1/2}-q^{1/2})q^{-s_i/2}=q^{-s_i},\qquad
 1-s_i(q^{-1/2}-q^{1/2})q^{s_i/2}=q^{s_i}.
\]
At the middle vertex for odd $N$, the two terms in parentheses
cancel. In the contraction sum of~\eqref{eq:R-ansatz}, substitution
of~\eqref{eq:c-closed} gives
\[
 (-1)^{p_j(p_i+p_j)}c_{ij}
 =(-1)^{p_ip_j}\vartheta_{i'}\vartheta_{j'}
                              q^{(\rho,\zeta_{i'}-\zeta_{j'})},
\]
using $\vartheta_{i'}=s_i\vartheta_i$ and $\zeta_{i'}=-\zeta_i$.
Reindexing by $i',j'$ and using
$(-1)^{p_ip_j}=s_j(-1)^{p_j(p_i+p_j)}$ combines the contraction and
exchange sums into the second line of~\eqref{eq:R}.
\end{proof}

Formula~\eqref{eq:R} is $R_0$ of \cite[(4.13)]{HT}, under the identification
\[
 \bar i_{\mathrm{HT}}=p_i,\qquad
 \varepsilon_i^{\mathrm{HT}}=\zeta_i,\qquad
 E_{ij}^{\mathrm{HT}}=e_{ij},\qquad s_{\mathrm{HT}}=n.
\]
Their parameters $m,n$ are $\dim V_{\bar0},\dim V_{\bar1}$,
respectively, and their parameter $q$ is our $q$.
The signs are those of~\eqref{eq:theta}.
Thus the Hasse recursion yields the same matrix in the same basis,
representation, and coproduct conventions.

Let $S=P R$ denote the braid matrix (recall that $P\in \End(V\tp V)$ is the super-permutation).
It is an invariant operator on $V\tp V$.
\begin{propn}\label{prop:inverse}
Under the substitution $q\mapsto\bar q$, we have
\begin{equation}\label{eq:inverse}
 R(q)^{-1}=R(\bar q).
\end{equation}
In particular,  $S^{-1}=P R_\infty$, where
\begin{equation}\label{eq:R-infinity}
 \begin{split}
 R_\infty={}&I+(q^{1/2}-q^{-1/2})\sum_i
       s_ie_{ii}\otimes
       \bigl(q^{(\zeta_i,\zeta_i)/2}e_{ii}
                      -q^{-(\zeta_i,\zeta_i)/2}e_{i'i'}\bigr)\\
 &+(q-\bar q)\sum_{i>j}s_je_{ij}\otimes
       \bigl(e_{ji}-(-1)^{p_j(p_i+p_j)}
          \vartheta_i\vartheta_jq^{(\rho,\zeta_i-\zeta_j)}e_{i'j'}\bigr).
 \end{split}
\end{equation}
\end{propn}
\begin{proof}
The matrices $\pi(e_i),\pi(f_i)$ are unchanged by $q\mapsto\bar q$,
including the terminal matrix of type $C$. Inverting the Cartan factors
interchanges $\Delta(x)$ and $\Delta^{\op}(x)$ for $x=e_i,f_i$. Hence
\[
 R(\bar q)\Delta^{\op}(x)=\Delta(x)R(\bar q),
\]
for these generators. For the Cartan generators, the same identity
follows from weight preservation. Thus $R(\bar q)R(q)$ commutes
with $\Delta(x)$ for every $x\in U_q(\g)$.
By~\eqref{eq:first-column}, its first column in the expansion over
the first tensor factor is $\delta_{i1}\id_V$.

The uniqueness argument of Proposition~\ref{prop:unique} also
applies to the commutation equation: comparison at $e_{ia}$
determines column $b$ from column $a$, since the coefficient
$\pi(f_k)_{ba}$ and the right Cartan factor $\pi(q^{-h_k/2})$
remain invertible. The identity operator satisfies the same equation
with the same initial column, so $R(\bar q)R(q)=I$.
Finally,
$S^{-1}=R(\bar q)P=P\bigl(P R(\bar q)P\bigr)$.
Interchanging the tensor factors in~\eqref{eq:R} and reindexing by
$i',j'$ gives~\eqref{eq:R-infinity}, which is also \cite[(4.14)]{HT}.
\end{proof}

On $V^{\otimes3}$, set
\[
 R_{12}=R\otimes\id_V,\qquad R_{23}=\id_V\otimes R,\qquad
 R_{13}=(\id_V\otimes P)(R\otimes\id_V)(\id_V\otimes P).
\]
Since Theorem~\ref{thm:R} computes the image of the normalized
universal $R$-matrix, applying $\pi^{\otimes3}$ to its
Yang--Baxter identity gives
\begin{equation}\label{eq:YBE}
 R_{12}R_{13}R_{23}=R_{23}R_{13}R_{12}.
\end{equation}
Moving the graded permutations $P$ through the even operators $R$
 gives the braid relation
$$S_{12}S_{23}S_{12}=S_{23}S_{12}S_{23}.$$

We next determine a polynomial relation for $S$, including the case
of coincident eigenvalues. Set $\delta=\sdim V=\sum_i s_i$.

\begin{thm}\label{thm:cubic}
The operator $S=P R$ satisfies
\begin{equation}\label{eq:cubic}
 (S-\bar q I)(S+qI)(S-q^{\delta-1}I)=0.
\end{equation}
If $\delta=0$, its minimal polynomial over $\C(q^{1/2})$ is
$(x-\bar q)^2(x+q)$ and $S$ is not semisimple.
\end{thm}
\begin{proof}
Consider the contraction operator
\[
 E=I+\frac{S-S^{-1}}{q-\bar q}.
\]
Formulas~\eqref{eq:R} and~\eqref{eq:R-infinity} give
\begin{equation}\label{eq:rank-one}
 E=\sum_{i,j}(-1)^{p_ip_j}\vartheta_i\vartheta_j
       q^{-(\rho,\zeta_i+\zeta_j)}e_{ij}\otimes e_{i'j'}.
\end{equation}
To verify this identity, use
\[
 P(e_{ab}\otimes e_{ba})=s_be_{bb}\otimes e_{aa},\qquad
 P(e_{ab}\otimes e_{a'b'})=s_ae_{a'b}\otimes e_{ab'}.
\]
The difference of the exchange terms, after multiplication by
$P$ and division by $q-\bar q$, is
$-\sum_{i\ne j}e_{ii}\otimes e_{jj}$. The diagonal terms give
\[
 -\sum_{i\ne i'}e_{ii}\otimes e_{ii}
                  +\sum_{i\ne i'}e_{i'i}\otimes e_{ii'},
\]
using $q^{s_i}-q^{-s_i}=s_i(q-\bar q)$.
Adding $I$ leaves the second sum and, for odd $N$, the middle term
$e_{n+1,n+1}\otimes e_{n+1,n+1}$. These are precisely the terms
of~\eqref{eq:rank-one} with $i=j'$.
For the contraction terms, set $i=a'$, $j=b$. Their coefficients are
\[
 s_a(-1)^{p_ap_b}\vartheta_a\vartheta_b
                          q^{(\rho,\zeta_a-\zeta_b)}
 =(-1)^{p_ip_j}\vartheta_i\vartheta_j
                          q^{-(\rho,\zeta_i+\zeta_j)},
\]
since $s_a\vartheta_a=\vartheta_{a'}$. The two triangular sums
cover all $a\ne b$, proving~\eqref{eq:rank-one}.

Taking the grading into account, the nonzero numerical block of $E$
has entries
\[
 s_j\vartheta_i\vartheta_jq^{-(\rho,\zeta_i+\zeta_j)}
 =\bigl(\vartheta_iq^{-(\rho,\zeta_i)}\bigr)
       \bigl(s_j\vartheta_jq^{-(\rho,\zeta_j)}\bigr).
\]
This factorization into a nonzero column and row shows that
$\operatorname{rank}E=1$.
As a Laurent polynomial in $S$, the operator $E$ commutes with $S$;
its one-dimensional image is therefore $S$-invariant, and $SE$
is proportional to $E$.

To find the scalar multiplier, compare coefficients of $e_{11}\otimes e_{1'1'}$.
In $E$, the coefficient is $s_1q^{-2(\rho,\zeta_1)}$.
For $SE=P RE$, we need the coefficient of
$e_{1'1}\otimes e_{11'}$ in $RE$.
Upper triangularity in the first tensor factor of $R$ leaves only
its diagonal term $q^{s_1}e_{1'1'}\otimes e_{11}$.
The coefficient of $e_{1'1}\otimes e_{11'}$ in $E$ is $1$,
and multiplication by $P$ contributes $s_1$.
Thus the coefficient in $SE$ is $s_1q^{s_1}$, whence
\[
 SE=q^{2(\rho,\zeta_1)+s_1}E=q^{\delta-1}E,
\]
where the last equality follows from~\eqref{eq:weyl}.
Multiplying the definition of $E$ by $(q-\bar q)S$ gives
\[
 (S-\bar q I)(S+qI)
 =S^2+(q-\bar q)S-I=(q-\bar q)q^{\delta-1}E.
\]
Multiplication by $S-q^{\delta-1}I$ proves~\eqref{eq:cubic}.

If $\delta=0$, then $S$ acts on the image of $E$ by $\bar q$.
Consequently,
\[
 E^2=\left(1+\frac{\bar q-q}{q-\bar q}\right)E=0,\qquad
 (S-\bar q I)(S+qI)=(q-\bar q)\bar q E\ne0.
\]
For distinct nonopposite indices $a<b$, the corresponding
two-dimensional weight block of $S$ is
\[
 s_ae_{ba}\otimes e_{ab}+s_be_{ab}\otimes e_{ba}
                         -(q-\bar q)e_{bb}\otimes e_{aa}.
\]
Its trace equals $-(q-\bar q)$ and its determinant is $-1$, so its
characteristic polynomial is
$x^2+(q-\bar q)x-1=(x-\bar q)(x+q)$.
Such a pair exists for $N>2$, and hence both eigenvalues occur.
The cubic polynomial in~\eqref{eq:cubic} is now
$(x-\bar q)^2(x+q)$; the nonzero quadratic expression above rules
out $(x-\bar q)(x+q)$. The cubic is therefore the minimal polynomial.
\end{proof}

\vspace{20pt}

\noindent
\underline{\large \bf Acknowledgement}

\vspace{10pt}

\noindent
This work is done at the Centre of Pure Mathematics, MIPT.

\noindent
D. Algethami is  thankful to the Deanship of Graduate Studies and Scientific Research at University of Bisha for supporting this work through the Fast-Track Research Support Program.

\noindent
N. Kryazhevskikh and A. Mudrov are   financially supported by Russian Science Foundation grant 26-11-00115.

\end{document}